\documentclass[leqno, 12pt]{amsart}

\usepackage{amssymb}
\usepackage{amsmath}
\newtheorem{theorem}{Theorem}[section]
\newtheorem{corollary}[theorem]{Corollary}

\newtheorem{lemma}[theorem]{Lemma}

\newtheorem{proposition}[theorem]{Proposition}

\theoremstyle{definition}
\newtheorem{remark}[theorem]{Remark}
\newtheorem{example}[theorem]{Example}
\newtheorem{definition}[theorem]{Definition}
\newtheorem{condition}[theorem]{Condition}

\begin{document}

\title{Contact Pairs}

\author{Gianluca Bande$^{1}$}
\address{$^{1}$ Universit\`{a} di Cagliari\endgraf Via Ospedale 72\endgraf 09129
Cagliari\endgraf Italia}
\email{gbande@unica.it}
\author{Amine Hadjar$^{2}$}
\address{$^{2}$ Universit\'{e} de Haute Alsace\endgraf 4 rue des fr\`{e}res
Lumi\`{e}re\endgraf 68093 Mulhouse Cedex\endgraf France}
\email{A.Hadjar@uha.fr}



\maketitle

\begin{abstract}

We introduce a new geometric structure on differentiable manifolds. A
\textit{Contact} \textit{Pair }on a $2h+2k+2$-dimensional manifold $M$ is a pair $\left(\alpha,\eta\right) $
of Pfaffian forms of constant classes $2k+1$ and
$2h+1$, respectively, whose characteristic foliations are transverse and
complementary and such that $\alpha$ and $\eta$ restrict to contact forms on the leaves of the characteristic foliation of $\eta$ and $\alpha$, respectively. Further differential objects are associated to Contact Pairs: two commuting Reeb vector fields, Legendrian curves on $M$ and two Lie
brackets on the set of differentiable functions on $M$. We give a local model
and several existence theorems on nilpotent Lie groups, nilmanifolds, bundles
over the circle and principal torus bundles.

\end{abstract}

\section{Introduction}

The aim of this paper is to study some differential Pfaffian forms of constant
class. This notion was introduced by E. Cartan (cf. \cite{CA2}, \cite{G1}). Global
problems relative to constant class forms have been efficiently studied in the
case of maximal class: contact forms, symplectic forms and generalized contact
forms (cf. \cite{B}).

In fact, we introduce a new geometric structure called \textit{Contact Pair}.
More precisely, a \textit{Contact Pair} (C.P.) of type $(h,k)$ on a $\left(
2h+2k+2\right)  $-dimensional manifold $M$ is a pair of Pfaffian forms
$(\alpha,\eta)$ satisfying the following properties:%

\begin{align*}
d\alpha^{h+1}  & =0\quad,\quad d\eta^{k+1}=0\quad\text{and}\\
& \alpha\wedge d\alpha^{h}\wedge\eta\wedge d\eta^{k}\text{ \quad is a volume
form on }M.
\end{align*}

The forms $\alpha$, $\eta$ have constant classes $2h+1$ and $2k+1$, respectively.

Some differential objects can be naturally associated to such a structure. The
\textit{characteristic foliations} of $\alpha$ and $\eta$ are transverse and
complementary. Their leaves are \textit{contact manifolds} of dimension $2k+1$
and $2h+1$, respectively. We give more general notions of a \textit{Reeb
vector field} and \textit{Legendrian curves}. We can also associate
\textit{two Lie brackets} on the algebra $\mathcal{C}^{\infty}\left(
M\right)  $. We show that contact pairs of the same type $(h,k)$ admit a local
model, like contact and symplectic forms.

Given the richness of this geometry, we are interested in the existence of
Contact Pairs. We give several existence theorems for nilpotent Lie groups,
nilmanifolds, bundles over the circle, as well as principal torus bundles
which showed their utility in contact geometry (\cite{L1}, \cite{L2},
\cite{H1}, \cite{H2}, \cite{B}). In the bundles $\left(  M^{4},B_{2}%
,T^{2}\right)  $, where the total space and the base are closed orientable
manifolds of dimensions 4 and 2, we construct $T^{2}$-invariant contact pairs
of type $(1,0)$.\smallskip

All geometric objects in this paper are supposed $C^{\infty}$.\smallskip

A similar structure called \textit{Contact-Symplectic Pair }is developed in
\cite{tesi} and \cite{artic 2 bande}.

\medskip

The authors would like to express
their gratitude to Y. Eliashberg and N. A' Campo for their
interest in this work, to M. Bordemann, M. Goze, R. Lutz and M.
Zessin for their valuable comments. They kindly acknowledge R.
Caddeo, S. Montaldo and P. Piu who made it possible for them to
meet.

\section{Contact Pairs (C.P.)}

Let $M$ be a $\left(  2h+2k+2\right)  $-dimensional manifold.

\begin{definition} \label{defCDC}
A pair $(\alpha$, $\eta)$ of Pfaffian forms on $M$ is said to be
a Contact Pair (for short C.P.) of type $(h,k)$ if the following conditions are
satisfied:
\begin{align*}
d\alpha^{h+1}  & =0\quad,\quad d\eta^{k+1}=0\\
& \text{and}\quad\alpha\wedge d\alpha^{h}\wedge\eta\wedge d\eta^{k}%
\quad\text{is a volume form on }M.
\end{align*}
\end{definition}

Thus the forms $\alpha$ and $\eta$ have constant classes $2h+1$ and $2k+1$,
respectively. A manifold equipped with a C.P. is clearly orientable.

A C.P. of type $(0,0)$ in a $2$-dimensional manifold $M$ is a pair of closed
Pfaffian forms with non-vanishing product; if $M$ is closed it follows that
$M$ is diffeomorphic to the 2-torus. Therefore we will always suppose $h\geq1
$ or $k\geq1$.

The simplest example of C.P.s is the following:

\bigskip

\textbf{\textquotedblleft Darboux\textquotedblright\ C.P.:} If $x_{1}$,$\cdots$, $x_{2h+1}$, $y_{1}$,$\cdots$, $y_{2k+1}$\ are coordinate functions on $\mathbb{R}^{2k+2h+2}$, then the forms
\[
\alpha=dx_{2h+1}+\sum_{i=1}^{h}x_{2i-1}dx_{2i}\quad,\quad\eta=dy_{2k+1}%
+\sum_{i=1}^{k}y_{2i-1}dy_{2i},
\]
(with the convention: if $h=0$ or $k=0$ the corresponding sum
is zero) determine a C.P. of type $(h,k)$ on $\mathbb{R}%
^{2k+2h+2}$.

\bigskip

This example is a local model of C.P.'s of type $(h,k)$ (see \S
\ref{parmodlocaCDC}).

\subsection{Reeb vector fields of a C.P}

In this section, we naturally generalize the notion of \textit{Reeb vector
field} classically associated to contact forms.

\begin{theorem}
\label{chdereebCDC}Let $(\alpha$, $\eta)$ be a \textup{C.P.} of type $(h,k)$ on $M$.
Then there exists a unique vector field $X_{\alpha}$ satisfying
\[
\alpha\left(  X_{\alpha}\right)  =1\quad,\quad i\left(  X_{\alpha}\right)
\left(  d\alpha^{h}\wedge\eta\wedge d\eta^{k}\right)  =0
\]
and a unique vector field $X_{\eta}$ satisfying
\[
\eta\left(  X_{\eta}\right)  =1\quad,\quad i\left(  X_{\eta}\right)  \left(
\alpha\wedge d\alpha^{h}\wedge d\eta^{k}\right)  =0\text{.}%
\]
\end{theorem}

\begin{proof}
For the uniqueness, suppose the existence of two vector fields $X_{\alpha}$
and $Y_{\alpha}$ verifying the first two equations. Then the volume form
$\alpha\wedge d\alpha^{h}\wedge\eta\wedge d\eta^{k}$ vanishes on $X_{\alpha
}-Y_{\alpha}$. Hence $X_{\alpha}=Y_{\alpha}$.

For the existence of $X_{\alpha}$, let us consider the form $\Omega
=d\alpha^{h}\wedge\eta\wedge d\eta^{k}$. Its characteristic space is
$1$-dimensional at every point, because $\Omega$\ is a non-vanishing $\left(
2k+2h+1\right)  $-form on a $\left(  2k+2h+2\right)  $-dimensional manifold.
Consider a tangent vector $u_{p}\neq0$ at a point $p$ such that $i\left(
u_{p}\right)  \Omega_{p}=0$. We set $\left(  X_{\alpha}\right)  _{p}%
=u_{p}/\alpha_{p}\left(  u_{p}\right)  $. This defines a smooth vector field
$X_{\alpha}$ on $M$ which satisfies the required conditions.
\end{proof}

\bigskip

For the \textit{Darboux} C.P., the Reeb vector fields are $X_{\alpha}%
=\partial/\partial x_{2h+1}$ and $X_{\eta}=\partial/\partial y_{2k+1}$\textit{.}

A simple computation shows the following additional properties of Reeb vector
fields of a C.P.

\begin{proposition}
The Reeb vector fields $X_{\alpha}$, $X_{\eta}$ of a \textup{C.P.} $(\alpha,\eta)$
commute and satisfy the following conditions:
\begin{align*}
\eta\left(  X_{\alpha}\right)   & =0,\quad i\left(  X_{\alpha}\right)
d\alpha=i\left(  X_{\alpha}\right)  d\eta=0,\\
\alpha\left(  X_{\eta}\right)   & =0,\quad i\left(  X_{\eta}\right)
d\alpha=i\left(  X_{\eta}\right)  d\eta=0.
\end{align*}
\end{proposition}

Hence $X_{\alpha}$\ (resp. $X_{\eta}$) is tangent to the characteristic
foliation of $\eta$\ (resp. $\alpha$), and coincides on every leaf with the
Reeb vector field (in the classical sense) of the contact form induced by
$\alpha$\ (resp. $\eta$) on the leaf.

\begin{corollary}
A \textup{C.P.} is invariant by the flows of its Reeb fields.
\end{corollary}

The following theorem shows that Reeb fields of C.P. have properties similar to Reeb fields of contact forms (cf. \cite{R}):

\begin{theorem}
The Reeb field $X_{\alpha}$ (resp. $X_{\eta}$) of a \textup{C.P.} $(\alpha,\eta)$ of
type $(h,k)$ on $M$ with $h\geqslant1$ (resp. $k\geqslant1$) does not admit
any closed transverse hypersurface.
\end{theorem}

\begin{proof}
As $h\geqslant1$, we have $d\alpha^{h}\wedge\eta\wedge d\eta^{k}=d\left(
\alpha\wedge d\alpha^{h-1}\wedge\eta\wedge d\eta^{k}\right)  $. If there
exists a closed transverse hypersurface $S$, then the form
\[
i\left(  X_{\alpha}\right)  \left(  \alpha\wedge d\alpha^{h}\wedge\eta
\wedge\eta^{k}\right)  =d\left(  \alpha\wedge d\alpha^{h-1}\wedge\eta\wedge
d\eta^{k}\right)
\]
is an exact volume form on $S$, which is impossible by Stokes' Theorem.
\end{proof}

\bigskip

The forms $d\alpha$\ and $d\eta$\ are absolute integral invariants of
$X_{\alpha}$\ and $X_{\eta}$. The forms $\alpha$\ and $\eta$\ are relative
integral invariants (cf. \cite{CA2}).

\subsection{Examples of C.P.'s}

\begin{enumerate}
\item Let $\left(  M_{1}^{2h+1}\text{, }\alpha\right)  $ and $\left(
M_{2}^{2k+1}\text{, }\eta\right)  $ be two contact manifolds and
$M=M_{1}^{2h+1}\times M_{2}^{2k+1}$. The pair $(\alpha$, $\eta)$ is a C.P. of
type $(h,k)$ on $M$ and it will be called \textit{Product} C.P. Its Reeb
fields are those of the two contact forms considered as vector fields on $M$.

For example, let $\theta_{1}$,$\theta_{2}$,$\theta_{3}$,$\phi$ be coordinate
functions on $\mathbb{R}^{4}$. Then the pair $\alpha=\sin\theta_{3}d\theta
_{1}-\cos\theta_{3}d\theta_{2}$, $\eta=d\phi$ is a C.P. of type $(1,0)$ on the
torus $T^{4}$ and its Reeb fields are $X_{\alpha}=\sin\theta_{3}%
\partial/\partial\theta_{1}-\cos\theta_{3}\partial/\partial\theta_{2}$,
$X_{\eta}=\partial/\partial\phi$.

\item Let $M_{1}^{2h+2}$ be a manifold with a C.P. $\left(  \alpha
,\eta\right)  $ of type $(h,0)$ and $M_{2}^{2k}$ an open manifold with a
volume form $\left(  d\theta\right)  ^{k}$, where $\theta$ is a Pfaffian form.
The pair $(\alpha$, $\eta+\theta)$ is a C.P. of type $(h,k)$ on $M_{1}%
^{2h+2}\times M_{2}^{2k}$.
\end{enumerate}

\subsection{Characteristic foliations of $\alpha$ and $\eta$}

Let $(\alpha$, $\eta)$ be a C.P. of type $(h,k)$ on $M$. We can naturally
associate to it the distribution of vectors on which $\alpha$ and $d\alpha$
vanish, and the one of vectors on which $\eta$ and $d\eta$ vanish. These
distributions are involutive because\emph{\ }$\alpha$ and $\eta$ have constant
classes. They determine the \textit{characteristic} foliations of\emph{\ }%
$\alpha$ and $\eta$, noted $\mathcal{F}$ and $\mathcal{G}$, respectively.

The\emph{\ }foliations\emph{\ }$\mathcal{F}$ and $\mathcal{G}$ are of
codimension $2h+1$ and $2k+1$, respectively, and \textit{their leaves are
contact manifolds}. This justifies the name of the structure. Moreover,
$\mathcal{F}$ and $\mathcal{G}$ are\textit{\ transverse and complementary}.

\section{Local model of contact pairs\label{parmodlocaCDC}}

To construct a local model for C.P., we use the existence of
characteristic foliations just described. One can easily show the
following (see \cite{tesi}).

\begin{theorem}
Let $(\alpha$, $\eta)$ be a \textup{C.P.} of type $\left(  h,k\right)  $ on $M$, with
$h\geqslant1$. For every point $p$ of $M$, there exists an open neighborhood
$V$ of $p$ and a coordinate system on $V$ such that the pair $(\alpha$,
$\eta)$ can be written:
\[
\alpha_{V}=dx_{2h+1}+\sum_{i=1}^{h}x_{2i-1}dx_{2i}\quad,\quad\eta
_{V}=dy_{2k+1}+\sum_{i=1}^{k}y_{2i-1}dy_{2i}\text{.}%
\]
with the convention: $\eta_{V}=dy_{1}$ if $k=0$.
\end{theorem}

Thus every C.P. is locally a product C.P. The open set $V$ will be called
a \textit{Darboux neighborhood}.

\section{Further differential objects associated to a C.P.}

Let $(\alpha$, $\eta)$ be a C.P. of type $(h,k)$ on a manifold $M$,
$\mathcal{F}$ and $\mathcal{G}$ the characteristic foliations of $\alpha$ and
$\eta$, respectively. We can naturally associate to it the following
differential objects:

\subsection{Characteristic foliations of $d\alpha$ and $d\eta$}

Since $d\alpha$ and $d\eta$ have constant classes $2h$ and $2k$, they
determine two characteristic foliations $\mathcal{F}^{\prime}$ and
$\mathcal{G}^{\prime}$of codimension $2h$ and $2k$, respectively.

Each leaf of $\mathcal{F}^{\prime}$ (resp. $\mathcal{G}^{\prime}$) is a union
of leaves of $\mathcal{F}$ (resp. $\mathcal{G}$). Furthermore, it is clear
that the pair induced by $(\alpha,\eta)$ on a leaf $F$ of $\mathcal{F}%
^{\prime}$ (resp. $\mathcal{G}^{\prime}$) is a C.P. of type $(0,k)$ (resp.
$(h,0)$) on $F.$ These foliations also have the following interesting properties:

\begin{proposition}
\label{thfeuilcardalfa}Suppose that the characteristic foliation
$\mathcal{F}^{\prime}$ of $d\alpha$ \textup{(}respectively $\mathcal{G}^{\prime}$ of
$d\eta$\textup{)} has a closed leaf $F$. Then all the leaves of $\mathcal{F}$ \textup{(}resp.
$\mathcal{G}$\textup{)} lying in $F$ are diffeomorphic, and $F$ fibers over the circle.
\end{proposition}

\begin{proof}
The form $\alpha$ induces on $F$ a non-zero closed Pfaffian form $\alpha_{F}$.
Then if $F$ is closed, it fibers over the circle (cf. \cite{Tichler}) and the
characteristic leaves of $\alpha_{F}$ are diffeomorphic. But these leaves are
exactly those of $\mathcal{F}$ lying in $F$.
\end{proof}

\subsection{Lie brackets on $\mathcal{C}^{\infty}\left(  M\right)  $
associated to a C.P}

Using the contact forms induced on the leaves of $\mathcal{F}$ and
$\mathcal{G}$, the algebra $\mathcal{C}^{\infty}\left(  M\right)  $ can be
endowed with a pair of Lie brackets. Precisely, to every function $f$ on $M$,
we can associate two vector fields $X_{f,\alpha}$ and $X_{f,\eta}$ as follows:

On each leaf $G$ of $\mathcal{G}$, there exists a unique vector field $X$
tangent to $G$ such that:
\[
\alpha_{G}\left(  X\right)  =f_{|G}\quad\text{and}\quad\left(  L_{X}\alpha
_{G}\right)  \wedge\alpha_{G}=0,
\]
where $\alpha_{G}$ is the contact form induced by $\alpha$ on $G$ (cf. \cite{LB}).

The vector field we obtained on $M$ is well defined, smooth and will be noted
$X_{f,\alpha}$. In a similar way, we construct $X_{f,\eta}.$ Now we can
introduce the two Lie brackets.

\begin{definition}
The Lie bracket of $f$, $g\in\mathcal{C}^{\infty}\left(  M\right)  $ along
$\alpha$ is the function:
\[
\left\{  f,g\right\}  _{\alpha}=\alpha\left(  \left[  X_{f,\alpha}%
,X_{g,\alpha}\right]  \right)  \text{,}%
\]
and the Lie bracket along $\eta$ is the function
\[
\left\{  f,g\right\}  _{\eta}=\eta\left(  \left[  X_{f,\eta},X_{g,\eta
}\right]  \right)  \text{.}%
\]
\end{definition}

The usual properties of Lie brackets hold.

\subsection{Legendrian curves}

They are defined as follows:

\begin{definition}
A Legendrian curve of the C.P. $(\alpha$, $\eta)$ with respect to $\alpha$ is
a piecewise differentiable curve $\gamma_{\alpha}$ on $M$ such that
\[
\alpha\left(  \overset{.}{\gamma}_{\alpha}\right)  =0\quad\text{and}\quad
i\left(  \overset{.}{\gamma_{\alpha}}\right)  \left(  \eta\wedge d\eta
^{k}\right)  \neq0\ \text{everywhere.}
\]
\end{definition}

Similarly, we define a Legendrian curve with respect to $\eta$. The curves
must be tangent to $\mathcal{F}$ (resp. $\mathcal{G}$), but transverse to
$\mathcal{G}$ (resp. $\mathcal{F}$). They can join the points as in connected
contact manifolds (see \cite{B}):

\begin{proposition}
Any two points on a connected manifold $M$ equipped with a \textup{C.P.} $(\alpha
,\eta)$ can be joined by a Legendrian curve with respect to $\alpha$ and by a
Legendrian curve with respect to $\eta$.
\end{proposition}

\section{Topological obstructions}

Let $(\alpha,\eta)$ be a C.P. of type $\left(  h,k\right)  $ on a manifold $M
$.

\begin{proposition}
If $M$ is a closed manifold, then $H^{2h+1}\left(  M,\mathbb{R}\right)  \neq0$
and $H^{2k+1}\left(  M,\mathbb{R}\right)  \neq0$.
\end{proposition}

\begin{proof}
If $H^{2h+1}\left(  M,\mathbb{R}\right)  =0$ or $H^{2k+1}\left(
M,\mathbb{R}\right)  =0$, the volume form $\alpha\wedge d\alpha^{h}\wedge
\eta\wedge d\eta^{k}$ is exact, which is impossible when $M$ is closed.
\end{proof}

\bigskip

An immediate consequence (which also follows from the existence of
a non-vanishing vector field on a C.P. manifold) is the following:

\begin{corollary}
There is no \textup{C.P.} on even-dimensional spheres.
\end{corollary}

By using \cite{Tichler}, we have the following result:

\begin{proposition}
\label{CPtypeh0 fibre}If $M$ is closed and equipped with a \textup{C.P.} of type
$(h,0)$, then $M$ fibers over the circle.
\end{proposition}

Here are some properties concerning Reeb vector fields:

\begin{proposition}
The Reeb vector fields of a \textup{C.P.} determine a locally free action
of $\mathbb{R}^{2}$. Every closed orbit is a $2$-torus.
\end{proposition}

\begin{proof}
As Reeb vector fields commute, they generate a locally free action
of $\mathbb{R}^{2}$. Every orbit admits two non-vanishing (Reeb)
fields and then its Euler-Poincar\'{e} characteristic vanishes.
\end{proof}

\smallskip

This action will be called the \textit{Reeb action}.

\begin{theorem}
The Reeb action does not admit a closed transversal submanifold of codimension
$2$.
\end{theorem}

\begin{proof}
If $h=k=0,$ it is obvious. Suppose $h\geq1$ and let $N$ be a $2$-codimensional
closed transversal submanifold. Then
\[
i\left(  Y\right)  i\left(  X\right)  \left(  \alpha\wedge d\alpha^{h}%
\wedge\eta\wedge d\eta^{k}\right)  =d\alpha^{h}\wedge d\eta^{k}=d\left(
\alpha\wedge d\alpha^{h-1}\wedge d\eta^{k}\right)
\]
induces an exact volume form on $N$, which is impossible by Stokes' Theorem.
\end{proof}

\begin{remark}
If every orbit of the Reeb action is a closed manifold, then we
have a locally free action of the torus. If this action generates
a principal fibre bundle, the C.P. (which is invariant) has an
empty singular set (see \S\ref{perCDCinvdim4} for
details).
\end{remark}

\section{\label{pargrdeLieCDC}C.P.'s on nilpotent Lie groups and nilmanifolds}

Nilpotent Lie groups and nilmanifolds provide further interesting examples of
C.P.'s. Below we present some examples of constructions in dimensions 4 and 6
which can possibly be extended to higher dimensions. We use the classification
of nilpotent Lie algebras of dimensions 4 and 6 in \cite{goze}.

\subsection{C.P.'s on nilpotent Lie groups}

In order to describe the Lie algebra of a Lie group, we give only the non-zero
ordered brackets of the fundamental fields $X_{i}$. Their dual forms will be
noted $\omega_{i}$.

\begin{example}
\textup{Consider the $4$-dimensional Lie algebra $n_{4}^{1}$ given by
\[
\left[  X_{1},X_{4}\right]  =X_{3}\quad,\quad\left[  X_{1},X_{3}\right]
=X_{2}\text{.}
\]
The pair $\left(  \omega^{2},\omega^{4}\right)  $ determines a C.P. of type
$\left(  1,0\right)  $ on the corresponding Lie group.}
\end{example}

\begin{example}
\textup{On the $6$-dimensional Lie algebra $n_{6}^{12}$ given by
\[
\left[  X_{1},X_{6}\right]  =X_{5}\quad,\quad\left[  X_{1},X_{5}\right]
=X_{4}\quad,\quad\left[  X_{2},X_{3}\right]  =X_{4}\text{,}
\]
the pair $\left(  \omega^{4},\omega^{6}\right)  $ determines a C.P. of type
$\left(  2,0\right)  $ on the corresponding Lie group.}
\end{example}

\begin{example}
\textup{On the group corresponding to the $6$-dimensional Lie algebra $n_{6}^{13}$
given by
\[
\left[  X_{1},X_{6}\right]  =X_{5}\text{, }\left[  X_{1},X_{5}\right]
=X_{4}\text{, }\left[  X_{1},X_{4}\right]  =X_{3}\text{, }\left[  X_{5}
,X_{6}\right]  =X_{2}\text{,}
\]
the pair $\left(  \omega^{2},\omega^{3}\right)  $ determines a C.P. of type
$\left(  1,1\right)  $.}
\end{example}

\subsection{C.P.'s on nilmanifolds}

We remark that in the previous examples the Lie algebras are rational; thus
the unique connected and simply connected Lie groups corresponding to them admit
cocompact discontinuous subgroups (cf. \cite{goze}). Then the
\textit{nilmanifolds,} obtained as quotients by these subgroups, are
\textit{closed manifolds equipped with} C.P.'s of the same type.

\section{Existence theorems of C.P.'s of type $(h,0)$}

As we have seen in \S\ref{CPtypeh0 fibre}, a closed manifold equipped with a
C.P. of type $(h,0)$ fibers over the circle. By using Feldbau's theorem
(cf. \cite{feldbau}), we construct non-product C.P.'s $(\alpha,\eta)$ of type
$(h,0)$ on manifolds that fiber over the circle, in such a way that the
characteristic foliation of $\eta$ coincides with the bundle foliation. We
recall Feldbau's theorem:

\textit{Equivalence classes of differentiable fiber bundles over the circle with closed, connected fiber $M$ and structural group $\mathrm{Diff}^{+}(M)$ are in one-to-one correspondence with $\pi_{0}(\mathrm{Diff}^{+}(M)).$}

If $f\in \mathrm{Diff}^{+}(M)$, the bundle is obtained as the quotient of
$M\times$ $]-\epsilon\,,1+\epsilon\lbrack$, $\epsilon>0$ by the
equivalence relation which identifies the points $(x,t)\in
M\times]-\epsilon\,,\epsilon\lbrack$ with $h(x,t)=(f(x),1+t)\in
M\times]1-\epsilon\,,1+\epsilon\lbrack$. The total space will be
denoted $M_{f}$ and $h$ will be called the \textit{gluing
diffeomorphism}.

\begin{theorem}
\label{thcontrfeldbau}Let $(B_{2h+1},\omega)$ be a connected, closed contact
manifold. If $f\in \mathrm{Diff}^{+}(B_{2h+1})$ and $f^{\ast}\omega=\omega$ , then
there exists a \textup{C.P.} $(\widetilde{\omega},\eta)$ of type $(h,0)$ on
$(B_{2h+1})_{f}$ . Moreover the pair can be chosen in such a way that every
contact leaf of the characteristic foliation of $\eta$ is a contact embedding
of $(B_{2h+1},\omega)$ .
\end{theorem}

\begin{proof}
Let $(B_{2h+1},\omega)$ be a connected, closed contact manifold,
$f\in \mathrm{Diff}^{+}(B)$ such that $f^{\ast}\omega=\omega$. Consider
$\epsilon>0$ and $I_{\epsilon}=]-\epsilon\,,1+\epsilon\lbrack$.
Let $p_{1},p_{2}$ be the projections of $B\times I_{\epsilon}$ on
$B$ and $I_{\epsilon}$, respectively, and $dt$ the canonical volume
form of $I_{\epsilon}$. The pair $(p_{1}^{\ast
}(\omega),p_{2}^{\ast}(dt))$ is a product C.P. of type $(h,0)$ on
$B\times I_{\epsilon}$, invariant by the \textit{gluing
diffeomorphism}. Thus it induces a C.P. $(\tilde{\omega},\eta)$ of
type $(h,0)$ on $B_{f}$. Let $\pi:B_{f}\rightarrow S^{1}$ be the
canonical projection, $d\theta$ the form on $S^{1}$ induced by
$dt$. By construction, we have $\eta=\pi^{\ast }d\theta$ and its
characteristic foliation coincides with the one defined by $\pi.$
Let $F=\pi^{-1}(\tau)$ be any fiber and $t\in I_{\epsilon}$ a
representative modulo 1 of $\tau.$ As a contactomorphism between
$\left( B,\omega\right)  $ and $\left(
F,\widetilde{\omega}_{F}\right)  $, one can take the one which
sends a point $p\in B$ to $\left(  p,t\right)  \in B\times
I_{\epsilon}$ modulo the gluing diffeomorphism $h$.
\end{proof}

\medskip

We shall say that the C.P. $(\tilde{\omega},\eta)$ constructed above is
\textit{induced by} $\omega$ and $f$. It is a product C.P. if and only if $f $
is isotopic to the identity map $\textup{id}_{B}$.

\begin{remark}
\label{metodeCP}This theorem gives a method to construct non-product C.P.'s on
a bundle over the circle where the fiber $B$ is endowed with a contact form
$\omega$ and a diffeomorphism $f$ leaving $\omega$ invariant. But not all
C.P.'s are obtained in this way, as the following example shows.
\end{remark}

Indeed, consider the forms
\[
\omega=\cos\theta_{3}d\theta_{1}+\sin\theta_{3}d\theta_{2}\quad\text{and}
\quad\eta=d\theta_{4}+\lambda d\theta_{1}
\]
on $T^{4}$, where $\lambda$ is an irrational number chosen
sufficiently small to ensure the pair $(\omega,\eta)$ to be a C.P.
and $\eta$ irrational. Therefore, the characteristic leaves of
$\eta$ are open. Thus they cannot be the compact fibers of a
bundle over the circle with total space $T^{4}$.

\begin{remark}
If the form $\eta$ is irrational it is close to a rational form
$\beta$. Then if $(\omega,\eta)$ is a C.P. of type $(h,0)$, so is
$(\omega,\beta)$. This shows that any C.P. of type $(h,0)$ is
close to a C.P. as in the previous theorem.
\end{remark}

\subsection{Non-product examples}

Here are some fundamental examples of pairs $\left(  \omega,f\right)  $ on a
manifold $B_{2h+1}$ which give C.P.'s on the associated Feldbau's bundle
$\left(  B_{2h+1}\right)  _{f}$.

\begin{enumerate}
\item On the torus $T^{3} $, for each integer $n\neq0 $, we consider the
contact form $\omega_{n} =\cos(n\theta_{3}) d\theta_{1} +\sin(n\theta_{3})
d\theta_{2}$ and the diffeomorphism $f_{n} (\theta_{1} ,\theta_{2} ,\theta
_{3})= (\theta_{2} ,\theta_{1} ,(\pi/2n) -\theta_{3})$.

\item We can also consider the contact form $\omega=\cos\theta_{3}d\theta
_{1}+\sin\theta_{3}d\theta_{2}$ and the diffeomorphism $f(\theta_{1}%
,\theta_{2},\theta_{3})=(\theta_{1},-\theta_{2},-\theta_{3})$ on $T^{3}$.

\item On the torus $T^{5}$, we have the contact form (see \cite{L2})
\begin{align*}
\omega & =\sin\theta_{2}\cos\theta_{2}d\theta_{1}-\sin\theta_{1}\cos\theta
_{1}d\theta_{2}+\\
& \cos\theta_{1}\cos\theta_{2}d\theta_{3}+(\sin\theta_{1}\cos\theta_{3}%
-\sin\theta_{2}\sin\theta_{3})d\theta_{4}+\\
& (\sin\theta_{1}\sin\theta_{3}+\sin\theta_{2}\cos\theta_{3})d\theta_{5}%
\end{align*}
and the diffeomorphism
\[
f(\theta_{1},\theta_{2},\theta_{3},\theta_{4},\theta_{5})=(\pi-\theta
_{1},-\theta_{2},\frac{\pi}{2}-\theta_{3},\theta_{5},\theta_{4})\text{.}%
\]

\item Let $U^{\ast}M$ be the unit cotangent bundle of an
$n$-dimensional Riemannian manifold $(M,g)$, $\alpha$ the
Liouville contact form on $U^{\ast }M$. If $f\in \mathrm{Diff}^{+}(M)$ has
finite order $p$, we choose the isomorphism $F$ of the bundle
$U^{\ast}M$ defined by
\[
F(x,\eta_{x})=(f(x),(f^{-1})^{\ast}(\eta_{x})/\left\Vert (f^{-1})^{\ast}%
(\eta_{x})\right\Vert _{g})
\]
for each $x\in M$ and $\eta_{x}\in U_{x}^{\ast}M.$ Clearly, its order is $p $
and $F^{\ast}\alpha=\lambda\alpha$ where $\lambda(x,\eta_{x})=1/\left\Vert
(f^{-1})^{\ast}(\eta_{x})\right\Vert _{g}$. Hence $\omega=\sum_{k=1}^{p}%
(F^{k})^{\ast}\alpha$ is a contact form on $U^{\ast}M$ which satisfies
$F^{\ast}\omega=\omega.$ The $\textup{C.P.}$ we obtain in $\left(  U^{\ast}M\right)
_{F}$ is of type $(n-1,0)$.
\end{enumerate}

The diffeomorphisms considered in (1), (2) and (3) do not induce the identity
map on the first homotopy group of the manifold. Hence, they are not isotopic
to $\textup{id}$. In (4), if $f$ is chosen non-isotopic to $\textup{id}$, so is $F $. In this
way we obtain non-product C.P.'s.

\subsection{Constructions on $(M_{3})_{f}$ where $M_{3}$ is the total space of
a principal $S^{1}$-bundle}

Let $B_{2}$ be a closed, connected, orientable surface of genus
$g\geq2.$ First, we construct a family $D(B_{2})$ of orientation preserving
diffeomorphisms of $B_{2}$ which are of finite order and not
isotopic to the identity map $\textup{id}_{B}$. Next, we intend to lift
these diffeomorphisms $f$ to certain principal $S^{1} $-bundles
$M_{3}$ over $B_{2},$ as isomorphisms $\widetilde{f}$ of the
bundle (also of finite order). Finally, we construct
$S^{1}$-invariant contact forms $\omega$ on $M_{3}$ satisfying
$\widetilde{f}^{\ast}\omega=\omega.$ Therefore,
we will have non-product C.P. of type $(1,0)$ on each bundle $(M_{3}%
)_{\widetilde{f}}$.

\subsubsection{\label{D(B2)}The family $D(B_{2})$}

Let $\Delta(g)=\{2\}\cup\{m\in\mathbb{N},m\geq2,m\mid(g-1)\}$ and
$n\in$ $\Delta(g)$ . We embed $B_{2}$ in $\mathbb{R}^{3}$ and give
diffeomorphisms $\varphi_{f,n}$ of order $n$ which are not
isotopic to $\textup{id}_{B}$, as follows:

\medskip

\noindent\textit{First case} ($n\mid(g-1$ $)$\textit{\ with}
$g\geq 3$)\textit{:} Let $l=(g-1)/n$. We consider a 2-torus $T$ of
revolution. Let $C_{1}$ be a meridian circle of $T$. By iterating
$n$ times the $(2\pi /n)$-rotation around the revolution axis, the
images of $C_{1}$ are $n$ new circles $C_{i}$ (except
$C_{n+1}=C_{1}$). We glue $l$ handles on one of the $n$ connected
components $T_{i}$ of $T-\bigcup_{i}C_{i}$ \ (after removing $2l $
disks). By each iteration of the same rotation, we will have $l$
new handles glued on $T$. Let $S^g _n$ be the surface so obtained
and $P_n =\{ f : B_2 \rightarrow S^g _n \mid f \; \textup{is a diffeomorphism} \}$. For each $f \in P_n$, this rotation induces a orientation preserving diffeomorphism $\varphi_{f, n}$ on $B_{2}$ of order $n$.

\medskip

\noindent\textit{Second case} ($n=2$\textit{\ with} $g$
\textit{even}): Consider the unit sphere $S$. Let $C_{1}$ be a
circle containing the poles, $l=g/n$, and $T_{i}$ the two
connected components of $S-C_{1}$. We glue $l$ handles on $T_{1}$
after removing $2l$ disks, and by the $\pi$-rotation sending
$T_{1}$ to $T_{2}$, we have $l$ new handles on $T_{2}$. As in the
first case, we obtain a surface $S^g _2$ and we put $P_2 =\{ f :
B_2 \rightarrow S^g _2 \mid f \; \textup{is a diffeomorphism} \}$.
The surface $S^g _2$ is also invariant by the above symmetry,
which therefore induces for every $f \in P_2$ a diffeomorphism
$\varphi_{f, 2}$ of order 2 on $B_{2}$. We put
$C_{2}=\varphi_{2}(C_{1})=C_{1}$ where $\varphi_{2}$ is the
$\pi$-rotation.

\medskip

By construction, for each $f \in P_n$ $B_{2}=f^{-1}
(\bigcup_{i}\Sigma_{i})$ where $\Sigma_{i}$ are $n$ 2-dimensional
compact connected submanifolds with boundary $C_{i}\cup C_{i+1} $,
and interior $T_{i}$ with the $l$ handles. For each $i$ and $f$,
$f^{-1} \circ \varphi _{f, n} \circ f (\Sigma_{i})=\Sigma_{i+1}$
and $\Sigma_{n+1}=\Sigma_{1}$.

We set $D(B_{2})=\{\varphi_{f, n}\mid n\in\Delta(g) , f \in
P_n\}$.

\subsubsection{How to lift these diffeomorphisms?}

In the trivial bundle $\left(  B_{2}\times S^{1},B_{2},S^{1}\right)  $, we can
evidently lift any diffeomorphism $f$ of $B_{2}$\ by $\widetilde{f}%
(x,\theta)=\left(  f(x),\theta\right)  $. However, it is possible
to find non-trivial principal $S^{1}$-bundles $\left(
M_{3},B_{2},S^{1}\right)  $ on which the elements of $D(B_{2})$
can be lifted as isomorphisms (of finite order) of the bundle. For
example, since $B_{2}$ was considered as a Riemannian submanifold
of $\mathbb{R}^{3},$ the elements $f$ of $D(B_{2})$ are isometries
(see \ref{D(B2)}). So, their tangent maps $\widetilde{f}$ induced
on the unit tangent bundle $UB_{2}$ are isomorphisms. We can also
consider for each integer $k\neq0$ the bundle $\left(
E_{(k)},B_{2},S^{1}\right)  $ associated to $\left(
UB_{2},B_{2},S^{1}\right)  $ with total space
$E_{(k)}=UB_{2}\times_{S^{1}}S^{1}$ quotient of $UB_{2}\times
S^{1}$ by the $S^{1}$-action $\left(  z,\left( p,z^{\prime}\right)\right)  \rightarrow \left( pz,z^{-k}z^{\prime}\right)$ (see \cite[p.54]{KN}). This bundle is still a principal
$S^{1}$-bundle. Every element $f\in D(B_{2})$ can be lifted as
follows: $\widehat{f_{k}}(\left[  p,z\right]  )=\left[ \widetilde
{f}(p),z\right]  $ for each representative $\left( p,z\right)  \in
UB_{2}\times S^{1}$ of an equivalence class $\left[  p,z\right]  $
in $E_{(k)}.$ For each $f\in D(B_{2}),$ the orientation preserving isomorphisms
$\widetilde{f}$ and $\widehat{f_{k}}$ have the same order as $f$,
and are not isotopic to the identity map.

\subsubsection{Construction of contact forms on $M_{3}$}

Let $M_{3}(B_{2},S^{1},q)$ be a principal $S^{1}$-bundle where the total space
$M_{3}$ is a closed, connected and orientable 3-manifold$.$ Let $f$ be any
element of $D(B_{2}).$ Assume that $F$ is a orientation preserving isomorphism of the bundle
with finite order inducing $f$ on the base. Then we have:

\begin{theorem}
There exists an $S^{1}$-invariant contact form $\omega$ on $M_{3}$ such that
$F^{\ast}\omega=\omega$.
\end{theorem}

\begin{proof}
In the above notation, let $n\in$ $\Delta(g)$ and $f$ $=\varphi_{n}$.

\smallskip

\noindent\textit{First case }($n\mid(g-1)$\textit{\ with
}$g\geq3$)\textit{:} Let $\eta_{1}$ be a germ of an
$S^{1}$-invariant contact form along the torus $q^{-1}(C_{1})$.
Next, we put $\eta_{2}=(F^{-1})^{\ast}\eta_{1}$. Then we have a
germ of an $S^{1}$-invariant contact form along $\partial\left(
q^{-1}(\Sigma_{1})\right)  $. According to \cite[Lemma 1.3 and Th. 3.3]{H1}, this extends
to an $S^{1}$-invariant contact form $\omega_{1}$ on
$q^{-1}(\Sigma_{1})$. Let $\omega$ be the form on $M_{3}$ whose
restriction to each $q^{-1}(\Sigma _{i})$ is
$\omega_{i}$=$(F^{-i+1})^{\ast}\omega_{1}$. It is well defined,
$S^{1}$-invariant and satisfies $F^{\ast}\omega=\omega$.

\smallskip

\noindent\textit{Second case }($n=2$\textit{\ with }$g$ \textit{even}): First,
let us construct a germ $\eta$ of an $S^{1}$-invariant contact form along the
torus $q^{-1}(C_{1})$ such that $F^{\ast}\eta=\eta$. There exists a tubular
neighborhood $U$ of this torus isomorphic to $V=T^{2}\times\left]
-1,1\right[  $, such that $F$ induces the isomorphism $F_{V}:(\theta
_{1},\theta_{2},t)\rightarrow(\theta_{1},\pi-\theta_{2},-t)$ of $\ V.$ The
$S^{1}$-invariant contact form $\widetilde{\eta}=d\theta_{1}+td\theta_{2}$ on
$V$ satisfies $F_{V}^{\ast}\widetilde{\eta}=\widetilde{\eta}.$ Hence we have
$\eta$. This germ also extends to an $S^{1}$-invariant contact form
$\omega_{1}$ on $q^{-1}(\Sigma_{1})$. Set $\omega_{2}=(F^{-1})^{\ast}%
\omega_{1}$ on $q^{-1}(\Sigma_{2})$, to get a global contact form on $M_{3}$
which coincides with $\omega_{i}$ on $q^{-1}(\Sigma_{i})$ for $i=1,2$.
\end{proof}
Thus $(M_{3})_{F}$ is endowed with a C.P. of type $(1,0)$.

\section{Existence of C.P.'s of type $(1,0)\ $on principal torus bundles
$\left(  M^{4},B_{2},T^{2}\right)  $\label{perCDCinvdim4}}

Several existence problems for contact forms have been solved using an
additional invariance condition under which the space carries geometrically
useful structures (see \cite{L1}, \cite{L2}, \cite{H1}, \cite{H2}, \cite{B}).
We proceed along the same lines. Let us consider principal torus bundles
$\left(  M^{4},B_{2},T^{2},\pi\right)  $ where the base and the total space
are connected, closed and orientable.

Let $\theta=\sum_{i=1}^{2}\theta^{i}\otimes e_{i}$ be a connection form,
$\Omega=\sum_{i=1}^{2}\Omega^{i}\otimes e_{i}$ its curvature form and $X_{1}
$, $X_{2}$ the fundamental vector fields generated by $e_{1}$, $e_{2}$.

Let $(\alpha$, $\eta)$ be a pair of Pfaffian forms on $M^{4}$. These forms are
$T^{2}$-invariant if and only if they can be written as follows:
\begin{align*}
\alpha & =\pi^{\ast}\left(  \beta\right)  +\pi^{\ast}\left(  f_{1}\right)
\theta^{1}+\pi^{\ast}\left(  f_{2}\right)  \theta^{2},\\
\eta & =\pi^{\ast}\left(  \gamma\right)  +\pi^{\ast}\left(  g_{1}\right)
\theta^{1}+\pi^{\ast}\left(  g_{2}\right)  \theta^{2},
\end{align*}
where $\beta$ and $\gamma$ are Pfaffian forms and $f_{i}$, $g_{i}$ functions
on the base space.

The pair $(\alpha$, $\eta)$ is an invariant C.P. of type $\left(  1,0\right)
$ if and only if the following three conditions\ are satisfied on $B_{2}$:

\begin{condition}
\label{cc1}
\begin{gather*}
\beta\wedge\left(  g_{2}df_{1}-g_{1}df_{2}\right)  +\left(  g_{2}f_{1}%
-g_{1}f_{2}\right)  \cdot\left(  d\beta+f_{1}\Omega^{1}+f_{2}\Omega^{2}\right)
\\
+\left(  f_{2}df_{1}-f_{1}df_{2}\right)  \wedge\gamma\neq0\quad\textup{
everywhere on } B_{2},
\end{gather*}
\end{condition}

\begin{condition}
\label{cc2}
\[
df_{1}\wedge df_{2}=0,
\]
\end{condition}

\begin{condition}
\label{cc3}
\begin{gather*}
dg_{1}=dg_{2}=0\quad\textup{(thus }g_{1}\textup{ and }g_{2}\textup{ are
constants),}\\
\textup{and\quad}d\gamma+g_{1}\Omega^{1}+g_{2}\Omega^{2}=0.
\end{gather*}
\end{condition}

We call the singular set of $(\alpha$, $\eta)$ the set $S$ on which the function
\[
\left(  \alpha\wedge\eta\right)  \left(  X_{1},X_{2}\right)  =\pi^{\ast
}\left(  g_{2}f_{1}-g_{1}f_{2}\right)  ,
\]
vanishes or equivalently $S=\pi^{-1}\left(  \Sigma\right)  $, where $\Sigma$
is the set of zeros of the function $h=g_{2}f_{1}-g_{1}f_{2}$. Geometrically,
$S$ is the set of points of $M^{4}$ where the forms $\alpha$ and $\eta$
induced on the fibers $T^{2}$ are proportional.

\subsection{Nature of the singular set}

If $\Sigma$ is not empty and does not coincide with $B_{2}$, one of the
constants $g_{1}$, $g_{2}$ is not zero. We suppose $g_{2}\neq0$, and by the 
condition \ref{cc1}, we have
\[
\left(  g_{2}\beta-f_{2}\gamma\right)  \wedge dh+g_{2}h\cdot\left(
d\beta+f_{1}\Omega^{1}+f_{2}\Omega^{2}-df_{2}\wedge\gamma\right)  \neq0.
\]

Where $h$ vanishes, we have $\left(  g_{2}\beta-f_{2}\gamma\right)
\wedge dh\neq0$, which implies that $\Sigma$ \textit{is a closed
orientable embedded submanifold of }$B_{2}$ \textit{of codimension
}$1$\textit{, }and therefore a finite disjoint union of circles.

The set $\Sigma$, when $\Sigma\neq\emptyset$ and $\Sigma\neq B_{2}$, verifies
the following obvious property:

\begin{condition}
\label{condenssing}\textup{To every connected component of $B_{2}-\Sigma
$ we can attach a sign + or - (the sign of the function $h$), in such a way that two adjacent components have opposite signs.}
\end{condition}

If $S=\pi^{-1}\left(  \Sigma\right)  $ is the singular set of an
invariant C.P., $\Sigma$ satisfies necessarily one of the three
following properties:

\begin{enumerate}
\item $\Sigma$ coincides with $B_{2}$,

\item $\Sigma$\ is empty,

\item $\Sigma$\textit{\ }is a $1$-codimensional closed orientable
embedded submanifold of $B_{2}$ satisfying the condition
\ref{condenssing}.
\end{enumerate}

We shall see in the next paragraph that each of these cases is possible.

\subsection{Existence theorems on $\left(  M^{4},B_{2},T^{2}\right)  $}

We shall show that for every set $\Sigma$ in each one of the three cases of
the previous paragraph there exists an invariant C.P. $(\alpha$, $\eta)$
having $\pi^{-1}(\Sigma)$ as singular set.

\subsubsection{The case $\Sigma=B_{2}$}

The base space is necessarily the 2-torus. In fact:

\begin{theorem}
Let $\left(  M^{4},B_{2},T^{2}\right)  $ be a principal torus bundle where
total and base space are closed, \textit{connected and orientable}. There
exists an invariant \textup{C.P.} of type $\left(  1,0\right)  $ with singular set
$M^{4}$ if and only if $B_{2}$ is the $2$-torus.
\end{theorem}

\begin{proof}
Suppose the existence of such a C.P. on $M^{4}$. Then we have
$\beta$, $\gamma$, $f_{i}$, $g_{i}$ on $B_{2}$ satisfying
the conditions \ref{cc1}, \ref{cc2}, \ref{cc3} and $h\equiv0$. These
conditions become

\begin{enumerate}
\item $\left(  f_{2}df_{1}-f_{1}df_{2}\right)  \wedge\gamma\neq0$,

\item $df_{1}\wedge df_{2}=0$,

\item $g_{1}=g_{2}=0$ and $d\gamma=0$.
\end{enumerate}

This implies that $\gamma$ is a non-singular closed Pfaffian form. Thus
$B_{2}$ is a torus.

Conversely, if the base space is a torus with \textquotedblleft pseudo-coordinates\textquotedblright \, $\theta^{1}$ and $\theta^{2}$, we can choose $\gamma=d\theta^{1}$, $f_{1}=\sin\theta^{2}$, $f_{2}=\cos\theta^{2}$ to get an invariant C.P. with
singular set $M^{4}$.
\end{proof}

\subsubsection{The case $\Sigma$ empty}

There is a constraint on the bundle. In fact:

\begin{theorem}
Let $\left(  M^{4},B_{2},T^{2}\right)  $ be a principal torus bundle with
total and base space closed, \textit{connected and orientable}. There exists
an invariant \textup{C.P.} of type $\left(  1,0\right)  $ with empty singular set if
and only if the characteristic classes of the bundle do not vanish simultaneously.
\end{theorem}

Hence in the trivial bundle $\left(  T^{2}\times T^{2},T^{2},T^{2}\right)  $,
there is no such C.P.

\bigskip

\begin{proof}
Suppose that there exists an invariant C.P. $(\alpha,\beta)$ with empty
singular set. In the notation above, at least one of $g_{1}$, $g_{2}$ is
not zero. If $g_{2}\neq0$, we have $f_{1}=\left(  g_{1}f_{2}+h\right)  /g_{2}%
$. Condition \ref{cc1} becomes
\[
g_{2}h^{2}d\left(  \beta/h\right)  +g_{2}h\cdot\left(  f_{1}\Omega^{1}%
+f_{2}\Omega^{2}\right)  +\left(  f_{2}dh-hdf_{2}\right)  \wedge\gamma
\neq0\text{,}%
\]
which gives
\[
g_{2}h\cdot d\left(  \beta/h\right)  -h\cdot d\left(  \left(  f_{2}/h\right)
\cdot\gamma\right)  +g_{2}\cdot\left(  f_{1}\Omega^{1}+f_{2}\Omega^{2}\right)
+f_{2}d\gamma\neq0\text{.}%
\]
Since $d\gamma=-g_{1}\Omega^{1}-g_{2}\Omega^{2}$, we have
\[
d\left(  g_{2}\beta/h\right)  -d\left(  \left(  f_{2}/h\right)  \cdot
\gamma\right)  +\Omega^{1}\neq0\quad\text{everywhere.}%
\]
This condition implies that $\int\Omega^{1}\neq0$ and thus $\left[  \Omega
^{1}\right]  \neq0$.

Conversely, suppose that there is a non-vanishing characteristic class (for
example $\left[  \Omega^{1}\right]  \neq0$) and consider the pair
\begin{align*}
\alpha & =\pi^{\ast}\left(  \beta\right)  +\pi^{\ast}\left(  g_{2}\right)
\theta^{1}+\pi^{\ast}\left(  -g_{1}\right)  \theta^{2}\text{,}\\
\eta & =\pi^{\ast}\left(  \gamma\right)  +\pi^{\ast}\left(  g_{1}\right)
\theta^{1}+\pi^{\ast}\left(  g_{2}\right)  \theta^{2}\text{,}
\end{align*}
where $g_{1}$, $g_{2}$ are constants and $\beta$, $\gamma$ Pfaffian forms on
$B_{2}$. This determines a C.P. of type (1,0) with empty singular set if and
only if
\begin{align*}
d\gamma+g_{1}\Omega^{1}+g_{2}\Omega^{2}  & =0\quad\text{and}\\
\left(  g_{2}^{2}+g_{1}^{2}\right)  \cdot\left(  d\beta+g_{2}\Omega^{1}
-g_{1}\Omega^{2}\right)   & \neq0\quad\text{everywhere.}
\end{align*}
As $\int\Omega^{1}\neq0$ we can find $\left(  g_{1},g_{2}\right)  \neq\left(
0,0\right)  $ such that $g_{1}\int\Omega^{1}+g_{2}\int\Omega^{2}=0$. Thus
$l=g_{2}\int\Omega^{1}-g_{1}\int\Omega^{2}\neq0$, and there exists a form
$\gamma$ satisfying $d\gamma+g_{1}\Omega^{1}+g_{2}\Omega^{2}=0$. Now choose a
volume form $\Phi$ whose integral is $1$. Then we have $\int l\Phi=g_{2}%
\int\Omega^{1}-g_{1}\int\Omega^{2},$ which implies the existence
of some $\beta$ satisfying
$d\beta+g_{2}\Omega^{1}-g_{1}\Omega^{2}\neq0$ everywhere.
\end{proof}

\subsubsection{$\Sigma$ is a submanifold of codimension $1$}

There is no obstruction on the bundle and we have:

\begin{theorem}
\label{th cdc dim 4}Consider a principal torus bundle $\left(  M^{4}%
,B_{2},T^{2},\pi\right)  $ with closed, connected, orientable
total and base space. Let $\Sigma$ be a closed orientable
$1$-codimensional imbedded submanifold of $B_{2}$, satisfying
the condition \textup{\ref{condenssing}}. There exists an invariant \textup{C.P.}
$\left(  \alpha,\eta\right)  $ of type $\left(  1,0\right)  $ on
$M^{4}$ having $\pi^{-1}\left(  \Sigma\right)  $ as the singular set.
\end{theorem}

An example is given in Remark \ref{metodeCP} where the singular set is the
union of two $3$-tori.

Before proving the theorem, we give a technical lemma which will be useful in
the following constructions (see \cite[p.5]{L1} for a proof).

\begin{lemma}
\label{lemmetechnique}Le $B$ be a closed orientable
$2$-dimensional manifold and $\Sigma$ a $1$-codimensional compact
orientable embedded submanifold satisfying the condition
 \textup{\ref{condenssing}}. Then there exists a function $h$ and a Pfaffian
form $\beta$ on $B$ such that $\Sigma=h^{-1}\left(  0\right)  $
and the form $h\cdot d\beta+\beta\wedge dh$ is a volume form on
$B$.
\end{lemma}

\begin{proof}
\textbf{(of Theorem \ref{th cdc dim 4})} We have to find $\beta$, $\gamma$,
$f_{i}$, $g_{i}$ such that the conditions \ref{cc1}, \ref{cc2} and \ref{cc3} are satisfied.

\smallskip

\noindent\textit{First step}: By Lemma \ref{lemmetechnique},\ there exist $h $
and $\beta_{0}$ such that
\[
\beta_{0}\wedge dh+h\cdot d\beta_{0}\neq0\text{, }h^{-1}\left(  0\right)
=\Sigma\text{.}
\]

\noindent\textit{Second step: }We construct constants $g_{i}$ and
a form $\gamma$ satisfying the condition \ref{cc3}. If the
characteristic classes $[\Omega^{i}]$ vanish, then
$\Omega^{1}=d\gamma_{1}$, $\Omega^{2}=d\gamma_{2}$. We take
$\left(  g_{1},g_{2}\right) \neq\left(  0,0\right)  $ arbitrarily
and $\gamma=-(g_{1}\gamma_{1}+g_{2}\gamma_{2})$.

If there is a non-vanishing characteristic class, we can always
find $\left( g_{1},g_{2}\right)  \neq\left(  0,0\right)  $ such
that: $g_{1}\int\Omega ^{1}+g_{2}\int\Omega^{2}=0$. This implies
that there exists $\gamma$ satisfying
$d\gamma+g_{1}\Omega^{1}+g_{2}\Omega^{2}=0$.
\smallskip

\noindent\textit{Third step: }We construct the functions $f_{i}$.
We set:
\[
f=h/\left(  k_{1}g_{2}-k_{2}g_{1}\right)  \quad,\text{\quad}f_{1}=k_{1}%
f\quad,\text{\quad}f_{2}=k_{2}f,
\]
where $k_{i}$ are constants such that
$k_{1}g_{2}-k_{2}g_{1}\neq0$. Then the condition \ref{cc2} is
satisfied.

\smallskip

\noindent\textit{Last step:} We set $\beta=r\beta_{0}$ where $r$ is a real
number sufficiently large such that the condition \ref{cc1} is satisfied.
\end{proof}

\begin{remark}
By Remark \ref{metodeCP}, if in the previous three theorems $\eta$ is
irrational then the C.P. $\left(  \alpha,\eta\right)  $ cannot be obtained as
in Theorem \ref{thcontrfeldbau}.
\end{remark}

\end{document}